\documentclass[reqno,11pt]{amsart}
\usepackage{amscd,amssymb,verbatim}
\numberwithin{equation}{section}
\theoremstyle{plain}

\newcommand\alp{\alpha}         

\newcommand\gam{\gamma}         \newcommand\Gam{\Gamma}
\newcommand\del{\delta}         \newcommand\Del{\Delta}
\newcommand\eps{\varepsilon}
\newcommand\zet{\zeta}
\newcommand\tet{\theta}

\newcommand\lam{\lambda}                \newcommand\Lam{\Lambda}

\newcommand\ome{\omega}         \newcommand\Ome{\Omega}

\newcommand\calB{{\mathcal{B}}}
\newcommand\calC{{\mathcal{C}}}

\newcommand\calI{{\mathcal{I}}}

\newcommand\calR{{\mathcal{R}}}

\newcommand\RR{\mathbb{R}}

\newcommand\ZZ{\mathbb{Z}}

\newcommand\CC{\mathbb{C}}

 \newcommand\grb{{\mathfrak{b}}}

 \newcommand\gro{{\mathfrak{o}}}

\newcommand\nek{,\ldots,}
\newcommand\sdp{\times \hskip -0.3em {\raise 0.3ex
\hbox{$\scriptscriptstyle |$}}} 

\newcommand\Det{\operatorname{Det}}

\newcommand\MOD{\operatorname{mod}}

\newcommand\Ker{\operatorname{Ker}}

\newcommand\rank{\operatorname{rank}}

\newcommand\RE{\operatorname{Re}}

\newcommand\Tr{\operatorname{Tr}}

\newcommand\oeta{{\overline{\eta}}}

\newcommand\hatL{{\widehat{L}}}

\newcommand\hatOme{{\widehat{\Ome}}}

\newcommand\tilB{{\widetilde{B}}}

\newcommand\tiln{{\tilde{n}}}

\newcommand\tilDel{{\widetilde{\Delta}}}

\newcommand\tilGam{{\widetilde{\Gamma}}}

\renewcommand{\>}{\rangle}
\newcommand{\<}{\langle}

\theoremstyle{plain}
\newtheorem{Thm}[subsection]{Theorem}
\newtheorem{Cor}[subsection]{Corollary}
\newtheorem{Lem}[subsection]{Lemma}
\newtheorem{Prop}[subsection]{Proposition}
\newtheorem{Conjec}[subsection]{Conjecture}

\newtheorem{Def}[subsection]{Definition}

\theoremstyle{remark}

\newtheorem{Rem}[subsection]{Remark}

\def\TeXref#1{%
        \leavevmode\vadjust{\setbox0=\hbox{{\tt
                \  {\tiny \textrm #1}}}%
        \theight=\ht0
        \advance\theight by \lineskip
        \kern -\theight \vbox to
        \theight{\rightline{\rlap{\box0}}%
        \vss}%
        }}%
\newif\ifShowLabels
\ShowLabelstrue
\newdimen\theight
\def\TeXrefEq#1{%
        \leavevmode\vadjust{\setbox0=\hbox{{\tt
                \  {\tiny \textrm #1}}}%
        \theight=\ht1
        \advance\theight by \lineskip
        \kern -\theight \vbox to
        \theight{\rightline{\rlap{\box0}}%
        \vss}%
        }}%

\newcommand{\refs}[1]{Section ~\ref{S:#1}}
\newcommand{\refss}[1]{Subsection ~\ref{SS:#1}}

\newcommand{\reft}[1]{Theorem ~\ref{T:#1}}
\newcommand{\refl}[1]{Lemma ~\ref{L:#1}}

\newcommand{\refconj}[1]{Conjecture ~\ref{Conj:#1}}
\newcommand{\refd}[1]{Definition ~\ref{D:#1}}

\newcommand{\refe}[1]{\eqref{E:#1}}

\newenvironment{thm}[1]%
        { \begin{Thm} \label{T:#1}  \ifShowLabels \TeXref{T:#1} \fi }%
        { \end{Thm} }

\renewcommand{\th}[1]{\begin{thm}{#1}  }
\renewcommand{\eth}{\end{thm} }

\newenvironment{lemma}[1]%
        { \begin{Lem} \label{L:#1}  \ifShowLabels \TeXref{L:#1} \fi }%
        { \end{Lem} }

\newcommand{\lem}[1]{\begin{lemma}{#1} }
\newcommand{\elem}{\end{lemma}}

\newenvironment{propos}[1]%
        { \begin{Prop} \label{P:#1}  \ifShowLabels \TeXref{P:#1} \fi }%
        { \end{Prop} }

\newcommand{\prop}[1]{\begin{propos}{#1} }
\newcommand{\eprop}{\end{propos}}

\newenvironment{corol}[1]%
        { \begin{Cor} \label{C:#1}  \ifShowLabels \TeXref{C:#1} \fi }%
        { \end{Cor} }
\newcommand{\cor}[1]{\begin{corol}{#1}  }
\newcommand{\ecor}{\end{corol}}

\newenvironment{conjec}[1]%
        { \begin{Conjec} \label{Conj:#1}  \ifShowLabels \TeXref{C:#1} \fi }%
        { \end{Conjec} }
\newcommand{\conj}[1]{\begin{conjec}{#1}  }
\newcommand{\econj}{\end{conjec}}

\newenvironment{defeni}[1]%
        { \begin{Def} \label{D:#1}  \ifShowLabels \TeXref{D:#1} \fi }%
        { \end{Def} }
\newcommand{\defe}[1]{\begin{defeni}{#1}  }
\newcommand{\edefe}{\end{defeni}}

\newenvironment{remark}[1]%
        { \begin{Rem} \label{R:#1}  \ifShowLabels \TeXref{R:#1} \fi }%
        { \end{Rem} }
\newcommand{\rem}[1]{\begin{remark}{#1}}
\newcommand{\erem}{\end{remark}}

\newcommand{\eq}[1]%
        { \ifShowLabels \TeXrefEq{E:#1} \fi
           \begin{equation} \label{E:#1} }
\newcommand{\eeq}{\end{equation}}
\newcommand{\meq}[1]%
        { \ifShowLabels \TeXrefEq{E:#1} \fi
           \begin{multline} \label{E:#1} }
\newcommand{\emeq}{\end{multline}}

\newcommand{\prf}{ \begin{proof} }
\newcommand{\eprf}{ \end{proof} }
\newcommand{\Label}[1]{\label{#1}  \ifShowLabels \TeXref{#1} \fi }

\ShowLabelsfalse

\newcommand{\n}{\nabla}

\renewcommand{\b}{\bullet}
\newcommand{\pa}{\text{\( \partial\)}}
\newcommand{\odd}{{\operatorname{odd}}}
\newcommand{\even}{{\operatorname{even}}}

\newcommand{\Arg}{\operatorname{\mathbf{Arg}}}
\newcommand{\Ph}{\operatorname{\mathbf{Ph}}}
\newcommand{\Mon}{\operatorname{Mon}}

\newcommand{\Detgrtetnp}{\Det_{{\operatorname{gr},\tet}}}

\newcommand{\Flat}{\operatorname{Flat}}\newcommand{\Eul}{\operatorname{Eul}}

\newcommand{\B}{\calB}\newcommand{\tB}{\tilde{\calB}}

\newcommand{\rat}{\rho_{\operatorname{an}}}

\renewcommand{\hatOme}{\check{\Ome}}

\newcommand{\PD}{\operatorname{PD}}
\newcommand{\PDp}{\operatorname{PD}'}

\newcommand{\trivial}{\operatorname{trivial}}
\renewcommand{\tiln}{\tilde{\n}}
\newcommand{\gr}{\operatorname{gr}}
\newcommand{\C}{{\calC}}

\newcommand{\tauBH}{\tau^{\operatorname{BH}}}
\begin{document}
\begin{flushright}
\end{flushright}

\title[Refined analytic and Burghelea-Haller torsions]{Comparison of the refined analytic and the Burghelea-Haller torsions}
\author[Maxim Braverman]{Maxim Braverman$^\dag$}
\address{Department of Mathematics\\
        Northeastern University   \\
        Boston, MA 02115 \\
        USA
         }
\email{maximbraverman@neu.edu} \begin{flushright}\em\small{To Yves Colin de Verdi\`ere at the occasion of the conference in his honor, in
admiration}\end{flushright}
\author[Thomas Kappeler]{Thomas Kappeler$^\ddag$}
\address{Institut fur Mathematik\\
         Universitat Z\"urich\\
         Winterthurerstrasse 190\\
         CH-8057 Z\"urich\\
         Switzerland
         }
\email{tk@math.unizh.ch}
\thanks{${}^\dag$The first author would like to thank the Max Planck Institute for Mathematics in Bonn, where part of this work was
completed.\\\indent ${}^\ddag$Supported in part by the Swiss National Science foundation, the programme SPECT, and the European Community through the
FP6 Marie Curie RTN ENIGMA (MRTN-CT-2004-5652)}

\begin{abstract}
The refined analytic torsion associated to a flat vector bundle  over a closed odd-dimensional manifold canonically defines a quadratic form
$\tau$ on the determinant line of the cohomology. Both $\tau$ and the Burghelea-Haller torsion are refinements of the Ray-Singer torsion. We
show that whenever the Burghelea-Haller torsion is defined it is equal to $\pm\tau$. As an application we obtain new results about the
Burghelea-Haller torsion. In particular, we prove a weak version of the Burghelea-Haller conjecture relating their torsion with the square of
the Farber-Turaev combinatorial torsion.
\end{abstract}
\maketitle

\section{Introduction}\Label{S:introduction}

\subsection{The refined analytic torsion}\label{SS:RAT}
Let $M$ be a closed oriented manifold of odd dimension $d=2r-1$ and let $E$ be a complex vector bundle over $M$ endowed with a flat connection $\n$.
In a series of papers \cite{BrKappelerRAT,BrKappelerRATdetline,BrKappelerRATdetline_hol}, we defined and studied the non-zero element
\[
    \rat\ = \ \rat(\n) \ \in \ \Det\big(\, H^\b(M,E)\,\big)
\]
of the determinant line $\Det\big(H^\b(M,E)\big)$ of the cohomology $H^\b(M,E)$ of $M$ with coefficients in $E$. This element, called the {\em
refined analytic torsion}, can be viewed as an analogue of the refinement of the Reidemeister torsion due to Turaev
\cite{Turaev86,Turaev90,Turaev01} and, in a more general context, to Farber-Turaev \cite{FarberTuraev99,FarberTuraev00}. The refined analytic
torsion carries information about the Ray-Singer metric and about the $\eta$-invariant of the odd signature operator associated to $\n$ and a
Riemannian metric on $M$. In particular, if $\n$ is a hermitian connection, then the Ray-Singer norm of $\rat(\n)$ is equal to 1. One of the
main properties of the refined analytic torsion is that it depends holomorphically on $\n$. Using this property we computed the ratio between
the refined analytic torsion and the Farber-Turaev torsion, cf. Th.~14.5 of \cite{BrKappelerRAT} and Th.~5.11 of
\cite{BrKappelerRATdetline_hol}. This result extends the classical Cheeger-M\"uller theorem about the equality between the Ray-Singer and the
Reidemeister torsions \cite{RaSi1,Cheeger79,Muller78,Muller93,BisZh92}.

\subsection{Quadratic form associated with $\rat$}\label{SS:ratqf}
We define the quadratic form $\tau= \tau_\n$ on the determinant line $\Det\big(H^\b(M,E)\big)$ by setting
\eq{ratqf}
    \tau(\rat) \ = \ e^{-\ 2\pi i\,\big(\,\eta(\n)-\rank E\cdot\eta_{\trivial}\,\big)},
\end{equation}
where $\eta(\n)$ stands for the $\eta$-invariant of the restriction to the even forms of the odd signature operator, associated to the flat vector bundle $(E,\n)$ and
a Riemannian metric on $M$ (cf. \refd{oddsign}), and $\eta_{\trivial}$ is the $\eta$-invariant of the trivial line bundle over $M$.

Properties of $\rat$, such as its metric independence or its analyticity established in \cite{BrKappelerRAT,BrKappelerRATdetline,BrKappelerRATdetline_hol} easily
translate into corresponding properties of $\tau_\n$ -- see \refss{applications}.

\rem{RATmodified}
The difference $\eta(\n)-\rank E\cdot\eta_{\trivial}$ in \refe{ratqf} is called the $\rho$-invariant of $(E,\n)$ and its reduction modulo $\ZZ$ is independent of the
Riemannian metric.
\erem
In a subsequent work we show that $\tau_\n$ can be defined directly, without going through the construction of $\rat$.

\subsection{The Burghelea-Haller complex Ray-Singer torsion}\label{SS:BHtorsion}
On a different line of thoughts, Burghelea and Haller \cite{BurgheleaHaller_function,BurgheleaHaller_function2} have introduced a refinement of
the {\em square} of the Ray-Singer torsion for a closed manifold of  arbitrary dimension, {\em provided that the complex vector bundle $E$
admits a non-degenerate complex valued symmetric bilinear form $b$}. They defined a complex valued quadratic form
\eq{tau}
    \tauBH\ =\ \tauBH_{b,\n}
\end{equation}
on the determinant line $\Det\big(H^\b(M,E)\big)$, which depends holomorphically on the flat connection $\n$ and is closely related to (the square
of) the Ray-Singer torsion. Burghelea and Haller then defined a complex valued quadratic form, referred to as {\em complex Ray-Singer torsion}. In
the case of a closed manifold $M$ of odd dimension it is given by
\eq{BHmodified}
    \tauBH_{b,\alp,\n} \ := \ e^{-2\,\int_M\,\ome_{\n,b}\wedge\alp}\cdot \tauBH_{b,\n},
\end{equation}
where $\alp\in \Ome^{d-1}(M)$ is an arbitrary closed $(d-1)$-form and $\ome_{\n,b}\in \Ome^1(M)$ is the Kamber-Tondeur form, cf. \cite[\S2]{BurgheleaHaller_function2}
-- see the discussion at the end of Section~5 of \cite{BurgheleaHaller_function2} for the reasons to introduce this extra factor. Burghelea and Haller conjectured
that, for a suitable choice of $\alp$, the form $\tauBH_{b,\alp,\n}$ is roughly speaking equal to the square of the Farber-Turaev torsion, cf.
\cite[Conjecture~5.1]{BurgheleaHaller_function2} and \refconj{BH} below.

Note that $\tauBH$ seems not to be related to the $\eta$-invariant, whereas the refined analytic torsion is closely related to it. In fact, our study of $\rat$ leads
to new results about $\eta$, cf. \cite[Th.~14.10, 14.12]{BrKappelerRAT} and \cite[Prop.~6.2, Cor.~6.4]{BrKappelerRATdetline_hol}.

\subsection{The comparison theorem}\label{SS:comparison}
The main result of this paper is the following theorem establishing a relationship between the refined analytic torsion and the Burghelea-Haller
quadratic form.
\th{RAT-BH}
Suppose $M$ is a closed oriented manifold of odd dimension $d=2r-1$ and let $E$ be a complex vector bundle over $M$ endowed with a flat connection $\n$. Assume that
there exists a symmetric bilinear form $b$ on $E$ so that the quadratic form \refe{tau} on $\Det\big(H^\b(M,E)\big)$ is defined. Then $\tauBH_{b,\n}= \pm\tau_\n$,
i.e.,
\eq{RAT-BH}
    \tauBH_{b,\n}\big(\rat(\n)\big) \ = \ \pm\, e^{-\ 2\pi i\,\big(\,\eta(\n)-\rank E\cdot\eta_{\trivial}\,\big)}.
\end{equation}
\eth
The proof is given in \refs{prRAT-BH}.

\reft{RAT-BH} implies that for manifolds of odd dimension, the inconvenient assumption of the existence of a non-degenerate complex valued symmetric bilinear form $b$
for the definition of the Burghelea-Haller torsion can be avoided, by defining the quadratic form via the refined analytic torsion as in \refe{ratqf}.

The relation between $\rat$ and $\tau$ (and, hence, when there exists a quadratic form $b$, with $\tauBH$) takes an especially simple form, when the bundle $(E,\n)$
is acyclic, i.e., when $H^\b(M,E)=0$. Then the determinant line bundle $\Det\big(H^\b(M,E)\big)$ is canonically isomorphic to $\CC$ and both, $\tau$ and $\rat$, can
be viewed as non-zero complex numbers. We then obtain the following
\cor{RAT-BH}
If\/  $H^\b(M,E)=0$, then
\eq{RAT-BHacyclic}
    \tau_\n \ = \ \Big(\,\rat(\n)\cdot e^{\pi i\,(\,\eta(\n)-\rank E\cdot\eta_{\trivial})}\,\Big)^{-2}.
\end{equation}
\ecor
In general, $\tau_\n$ (and, hence, $\tau_\n^{\operatorname{BH}}$) does not admit a square root which is holomorphic in $\n$, cf. Remark~5.12 and
the discussion after it in \cite{BurgheleaHaller_function2}. In particular, the product\/ $\rat\cdot{}e^{\pi i(\eta(\n)-\rank
E\cdot\eta_{\trivial})}$ is not a holomorphic function of $\n$, since $e^{\pi i(\eta(\n)-\rank E\cdot\eta_{\trivial})}$ is not even continuous
in $\n$. Thus the refined analytic torsion can be viewed as a modified version of the inverse square root of $\tau_\n$, which is holomorphic.

\subsection{Properties of the quadratic forms $\tau$ and $\tauBH$}\label{SS:applications}
As an application of our previous papers \cite{BrKappelerRAT,BrKappelerRATdetline,BrKappelerRATdetline_hol} we obtain various results about the quadratic form $\tau$,
some of them generalizing known properties of the Burghelea-Haller torsion $\tauBH$.  In particular, we show that $\tau$ is independent of the choice of the
Riemannian metric. As an application of \reft{RAT-BH} one sees that $\tauBH_{b,\alp,\n}$ is invariant under the deformation of the non-degenerate bilinear form $b$
(cf. \reft{BHanomaly}) -- a result, which was first proven by Burghelea and Haller  \cite[Th.~4.2]{BurgheleaHaller_function2}. We also slightly improve this result,
cf. \reft{BHanomaly2}.

Next we discuss our main application of \reft{RAT-BH}.
\subsection{Comparison between the Farber-Turaev and the Burghelea-Haller torsions}\label{SS:BH-Tur}
In \cite{BurgheleaHaller_function2}, Burghelea and Haller made a conjecture relating the quadratic form \refe{BHmodified} with the refinement of
the combinatorial torsion introduced by Turaev \cite{Turaev86,Turaev90,Turaev01} and, in a more general context, by Farber and Turaev
\cite{FarberTuraev99,FarberTuraev00}, cf. \cite[Conjecture~5.1]{BurgheleaHaller_function2}. In the case when the bundle $E$ is acyclic and
admits a non-degenerate complex valued symmetric bilinear form, the Burghelea-Haller conjecture states that \refe{BHmodified} is equal to the
square of the Turaev torsion. More precisely, recall that the Turaev torsion depends on the Euler structure $\eps$ and a choice of a
cohomological orientation, i.e, an orientation $\gro$ of the determinant line of the cohomology $H^\b(M,\RR)$ of $M$. The set of Euler
structures $\Eul(M)$, introduced by Turaev, is an affine version of the integer homology $H_1(M,\ZZ)$ of $M$. It has several equivalent
descriptions \cite{Turaev86,Turaev90,Burghelea99,BurgheleaHaller_Euler}. For our purposes, it is convenient to adopt the definition from
Section~6 of \cite{Turaev90}, where an Euler structure is defined as an equivalence class of nowhere vanishing vector fields on $M$ -- see
\cite[\S5]{Turaev90} for the description of the equivalence relation. The definition of the Turaev torsion was reformulated by Farber and Turaev
\cite{FarberTuraev99,FarberTuraev00}. The Farber-Turaev torsion, depending on $\eps$, $\gro$, and $\n$, is an element of the determinant line
$\Det\big(H^\b(M,E)\big)$, which we denote by $\rho_{\eps,\gro}(\n)$.

Recall that the quadratic form $\tauBH_{b,\alp,\n}$ is defined in \refe{BHmodified}. Burghelea and Haller made a conjecture,
\cite[Conjecture~5.1]{BurgheleaHaller_function2}, relating the quadratic form $\tauBH_{b,\alp,\n}$ and $\rho_{\eps,\gro}(\n)$, which extends the Bismut-Zhang theorem
\cite{BisZh92}.  They have proved their conjecture modulo sign in the case when the dimension of the manifold $M$ is even and the bundle $E$ admits a parallel
symmetric bilinear form  (\cite[Th.~5.7]{BurgheleaHaller_function2}) and in some other cases as well - see \cite[\S5]{BurgheleaHaller_function2}. Though Burghelea and
Haller stated their conjecture for manifolds of arbitrary dimensions, we restrict our formulation to the odd dimensional case.

Suppose $M$ is a closed oriented odd dimensional manifold. Let $\eps\in \Eul(M)$ be an Euler structure on $M$ represented by a non-vanishing
vector field $X$. Fix a Riemannian metric $g^M$ on $M$ and let $\Psi(g^M)\in \Ome^{d-1}(TM\backslash\{0\})$ denote the Mathai-Quillen form,
\cite[\S7]{MathaiQuillen}, \cite[pp.~40-44]{BisZh92}. Set
\[
        \alp_\eps \ = \ \alp_\eps(g^M) \ := \ X^*\Psi(g^M)\ \in \ \Ome^{d-1}(M).
\]
This is a closed differential form, whose cohomology class $[\alp_\eps]\in H^{d-1}(M,\RR)$ is closely related to the integer cohomology class,
introduced by  Turaev \cite[\S5.3]{Turaev90} and called {\em the characteristic class $c(\eps)\in H_1(M,\ZZ)$ associated to an Euler structure
$\eps$}. More precisely, let $\PD:H_1(M,\ZZ)\to H^{d-1}(M,\ZZ)$ denote the Poincar\'e isomorphism. For $h\in H_1(M,\ZZ)$ we denote by $\PDp(h)$ the
image of $\PD(h)$ in $H^{d-1}(M,\RR)$. Then
\eq{IPDceps}
    \PD'\big(\,c([X])\,\big) \ = \ -2\,[\alp_\eps] \ = \  -\,2\,[X^*\Psi(g^M)],
\end{equation}
and, hence,
\eq{ome=[ome]}
    2\,\int_M\,\ome_{\n,b}\wedge\alp_\eps \ = \ -\,\langle\,[\ome_{\n,b}],c(\eps)\,\rangle,
\end{equation}
where $\ome_{\n,b}\in \Ome^1(M)$ is the Kamber-Tondeur form, cf. \refe{BHmodified}.

Note that \refe{IPDceps} implies that $2\alp_\eps$ represents an integer class in $H^{d-1}(M,\RR)$.

\conj{BH}{\rm\textbf{[Burghelea-Haller]}} Assume that $(E,\n)$ is a flat vector bundle over $M$ which admits a non-degenerate symmetric bilinear
form $b$. Then
\eq{BHconj}
    \tauBH_{b,\alp_\eps,\n}\big(\,\rho_{\eps,\gro}(\n)\,\big) \ = \ 1,
\end{equation}
or, equivalently,
\eq{BHconj1}
    \tauBH_{b,\n}\big(\,\rho_{\eps,\gro}(\n)\,\big) \ = \ e^{2\,\int_M\,\ome_{\n,b}\wedge\alp_\eps}.
\end{equation}
\econj

\subsection{A generalization of the Burghelea-Haller conjecture}\Label{SS:genBH}
Following Farber \cite{Farber00AT}, we denote by $\Arg_\n$ the unique cohomology class $\Arg_\n\in H^1(M,\CC/\ZZ)$ such that for every closed curve
$\gam$ in $M$ we have
\eq{IIArg}
    \det\big(\,\Mon_\n(\gam)\,\big) \ = \ \exp\big(\, 2\pi i\<\Arg_\n,[\gam]\>\,\big),
\end{equation}
where $\Mon_\n(\gam)$ denotes the monodromy of the flat connection $\n$ along the curve $\gam$ and $\<\cdot,\cdot\>$ denotes the natural pairing
$H^1(M,\CC/\ZZ)\,\times\, H_1(M,\ZZ) \to  \CC/\ZZ$.

By Lemma~2.2 of \cite{BurgheleaHaller_function2} we get
\eq{omeseps}
    e^{-\<[\ome_{\n,b}],c(\eps)\>} \ = \ \pm\, \det\,\Mon_\n\big(c(\eps)\big) \ = \ \pm\, e^{2\pi i\<\Arg_\n,c(\eps)\>}.
\end{equation}
(Note that $\Mon_\n(\gam)$ is equal to the {\em inverse} of what is denoted by $\operatorname{hol}^E_x(\gam)$ in \cite{BurgheleaHaller_function2}).

Combining \refe{ome=[ome]}, \refe{IIArg} and \refe{omeseps} we obtain
\[
     e^{2\,\int_M\,\ome_{\n,b}\wedge\alp_\eps} \ = \ \pm\, e^{2\pi i\<\Arg_\n,c(\eps)\>}.
\]
Thus, up to sign, the Burghelea-Haller conjecture \refe{BHconj1} can be rewritten as
\eq{BHconj2}
    \tauBH_{b,\n}\big(\,\rho_{\eps,\gro}(\n)\,\big) \ = \ \pm\, e^{2\pi i\<\Arg_\n,c(\eps)\>}.
\end{equation}
In view of \reft{RAT-BH} we make the following stronger conjecture involving $\tau_\n$ instead of $\tauBH_{b,\alp_\eps,\n}$, and, hence, meaningful also in the
situation, when the bundle $E$ does {\em not} admit a non-degenerate symmetric bilinear form.

\conj{BHnew} Assume that $(E,\n)$ is a flat vector bundle over $M$. Then
\eq{BHconjnew}
    \tau_\n\big(\,\rho_{\eps,\gro}(\n)\,\big) \ = \  e^{2\pi i\<\Arg_\n,c(\eps)\>},
\end{equation}
or, equivalently,
\eq{BHconjnew1}
    e^{ \pi i\,\big(\,\eta(\n)-\rank E\cdot\eta_{\trivial}\,\big)}\cdot\rat(\n) \ = \ \pm\, e^{-\pi i\<\Arg_\n,c(\eps)\>}\cdot \rho_{\eps,\gro}(\n).
\end{equation}
 \econj

\refconj{BHnew} implies the Burghelea-Haller conjecture~\ref{Conj:BH} up to sign.

\rem{BHconjnew1}
By construction, the left hand side of \refe{BHconjnew1} is independent of the Euler structure $\eps$ and the cohomological orientation $\gro$, while the right hand
side of \refe{BHconjnew1} is independent of the Riemannian metric $g^M$. Note that the fact that $e^{ \pi i\,\big(\,\eta(\n)-\rank
E\cdot\eta_{\trivial}\,\big)}\cdot\rat(\n)$ is independent of $g^M$ up to sign follows immediately from Lemma~9.2 of \cite{BrKappelerRATdetline}, while the fact that
$e^{-\pi i\<\Arg_\n,c(\eps)\>}\cdot \rho_{\eps,\gro}(\n)$ is independent of $\eps$ and independent of $\gro$ up to sign is explained on page~212 of
\cite{FarberTuraev00}.
\erem

In Theorem~5.1 of \cite{BrKappelerRATdetline_hol} we computed the ratio of the refined analytic and the Farber-Turaev torsions. Using this
result and \reft{RAT-BH} we establish the following weak version of \refconj{BHnew} (and, hence, of the Burghelea-Haller
conjecture~\ref{Conj:BH}).

\th{BHconj}
(i) \ Under the same assumptions as in \refconj{BHnew}, for each connected component $\C$ of the set\/ $\Flat(E)$ of\/ flat connections on $E$ there exists a constant
$R_\C$ with $|R_\C|=1$, such that
\eq{BHconjweak}
    \tau_\n\big(\,\rho_{\eps,\gro}(\n)\,\big) \ = \ R_\C\cdot e^{2\pi i\<\Arg_\n,c(\eps)\>}, \qquad\text{for all}\quad \n\in \C.
\end{equation}

(ii) \ If the connected component $\C$ contains an acyclic Hermitian connection then $R_\C= 1$, i.e.,
\eq{BHconjweak1}
    \tau_\n\big(\,\rho_{\eps,\gro}(\n)\,\big) \ = \  e^{2\pi i\<\Arg_\n,c(\eps)\>}, \qquad\text{for all}\quad \n\in \C.
\end{equation}
\eth
The proof is given in \refss{comparisonBHTUR}.

\rem{Haller}
(i) \ The second part of \reft{BHconj} is due to Rung-Tzung Huang, who also proved it in the case when $\C$ contains a Hermitian connection
which is not necessarily acyclic, \cite{Huang06}.

(ii) \ It was brought to our attention by Stefan Haller that one can modify the arguments of our proofs of \reft{RAT-BH} and of
\cite[Th.~5.1]{BrKappelerRATdetline_hol} so that they can be applied directly to the Burghelea-Haller torsion. It might lead to a proof of an analogue of
\reft{BHconj} for $\tauBH_{\n,b}$ on an even dimensional manifold.
\erem

\subsection{Added in proofs}\Label{SS:addedinproofs}
Since the first version of our paper was posted in the archive a lot of progress has been made. First, Huang \cite{Huang06} showed that if the connected component
$\calC\subset \Flat(E)$ contains a Hermitian connection, then the constant $R_\C$ of \reft{BHconj} is equal to $1$. Part of his result is now incorporated in item
(ii) of our \reft{BHconj}. Later Burghelea and Haller (D.~Burghelea and S.~Haller, \emph{{Complex valued Ray--Singer torsion II}}, \texttt{arXiv:math.DG/0610875})
proved \refconj{BH} up to sign. Independently and at the same time Su and Zhang (G.~Su and W.~Zhang, \emph{{A Cheeger-Mueller theorem for symmetric bilinear
torsions}}, \texttt{arXiv:{}math.DG/0610577}) proved \refconj{BH} in full generality. Both proofs used methods completely different from ours. In fact,
Burghelea-Haller, following \cite{BFK3}, and Su-Zhang, following \cite{BisZh92}, study a Witten-type deformation of the non-self adjoint Laplacian \refe{Lapl} and
adopt all arguments of these papers to the new situation. In contrast, our \reft{BHconj} provides a ``low-tech" approach to the Burghelea-Haller conjecture and, more
generally, to \refconj{BHnew}. On the other side, it would be interesting to see if the methods of Burghelea-Haller and Su-Zhang can be used to prove \refconj{BHnew}.

\subsection*{Acknowledgment} We would like to thank Rung-Tzung Huang for suggesting to us the second part of \reft{BHconj}. We are also grateful to Stefan Haller for valuable comments on a preliminary version of this paper.

\section{The Refined Analytic Torsion}\label{S:RAT}

In this section we recall the definition of the refined analytic torsion from \cite{BrKappelerRATdetline}. The refined analytic torsion is
constructed in 3 steps: first, we define the notion of refined torsion of a finite dimensional complex endowed with a chirality operator, cf.
\refd{refinedtorsion}. Then we  fix a Riemannian metric $g^M$ on $M$ and consider the odd signature operator $\B = \B(\n,g^M)$ associated to a flat
vector bundle $(E,\n)$, cf. \refd{oddsign}. Using the {\em graded determinant} of $\B$ and the definition of the refined torsion of a finite
dimensional complex with a chirality operator we construct an element $\rho= \rho(\n,g^M)$ in the determinant line of the cohomology, cf. \refe{rho}.
The element $\rho$ is almost the refined analytic torsion. However, it might depend on the Riemannian metric $g^M$ (though it does not if
$\dim{}M\equiv1\ (\MOD 4)$ or if $\rank(E)$ is divisible by 4). Finally we ``correct" $\rho$ by multiplying it by an explicit factor, the metric
anomaly of $\rho$, to obtain a diffeomorphism invariant $\rat(\n)$ of the triple $(M,E,\n)$, cf. \refd{rat}.

\subsection{The determinant line of a complex}\label{S:detline}
Given a complex vector space $V$ of dimension $\dim{}V= n$, the {\em determinant line} of $V$ is the line $\Det(V):= \Lam^nV$, where $\Lam^nV$
denotes the $n$-th exterior power of $V$. By definition, we set $\Det(0):= \CC$. Further, we denote by $\Det(V)^{-1}$ the dual line of $\Det(V)$. Let
\eq{Cpintrod}
    \begin{CD}
       (C^\b,\pa):\quad  0 \ \to C^0 @>{\pa}>> C^1 @>{\pa}>>\cdots @>{\pa}>> C^d \ \to \ 0
    \end{CD}
\end{equation}
be a  complex of finite dimensional complex vector spaces. We call the integer $d$ the {\em length} of the complex $(C^\b,\pa)$ and denote by
$H^\b(\pa)=\bigoplus_{i=0}^d H^{i}(\pa)$ the cohomology of $(C^{\b},\pa)$. Set
\eq{DetCH}
    \Det(C^\b) \ := \ \bigotimes_{j=0}^d\,\Det(C^j)^{(-1)^j}, \qquad \Det(H^\b(\pa))\ := \ \bigotimes_{j=0}^d\,\Det(H^j(\pa))^{(-1)^j}.
\end{equation}
The lines $\Det(C^\b)$ and $\Det(H^\b(\pa))$ are referred to as the {\em determinant line of the complex} $C^\b$ and the {\em determinant line of its
cohomology}, respectively. There is a canonical isomorphism
\eq{isomorphism}
     \phi_{C^\b} \ = \ \phi_{(C^\b,\pa)}:\, \Det(C^\b) \ \longrightarrow \ \Det(H^\b(\pa)),
\end{equation}
cf., for example, \S2.4 of \cite{BrKappelerRATdetline}.

\subsection{The refined torsion of a finite dimensional complex with a chirality operator}\label{SS:rtfd}
Let $d= 2r-1$ be an odd integer and let $(C^\b,\pa)$ be a length $d$ complex of finite dimensional complex vector spaces. A {\em chirality operator}
is an involution $\Gam:C^\b\to C^\b$ such that $\Gam(C^j)= C^{d-j}$, $j=0\nek d$. For $c_j\in \Det(C^j)$ $(j=0\nek d)$ we denote by $\Gam{}c_j\in
\Det(C^{d-j})$ the image of $c_j$ under the isomorphism $\Det(C^j)\to \Det(C^{d-j})$ induced by $\Gam$. Fix non-zero elements $c_j\in \Det(C^j)$,
$j=0\nek r-1$ and denote by $c_j^{-1}$ the unique element of $\Det(C^j)^{-1}$ such that $c^{-1}_j(c_j)=1$. Consider the element
\eq{Gamc}
    c_{{}_\Gam} \ := \ (-1)^{\calR(C^\b)}\cdot
    c_0\otimes c_1^{-1}\otimes \cdots \otimes c_{r-1}^{(-1)^{r-1}}\otimes (\Gam c_{r-1})^{(-1)^r}\otimes (\Gam c_{r-2})^{(-1)^{r-1}}
    \otimes \cdots\otimes (\Gam c_0)^{-1}
\end{equation}
of $\Det(C^\b)$, where

\eq{R(C)}
  \calR(C^\b) \ := \ \frac12\ \sum_{j=0}^{r-1}\, \dim C^{j}\cdot\big(\, \dim C^j+(-1)^{r+j}\,\big).
\end{equation}
It follows from the definition of $c_j^{-1}$ that $c_{{}_\Gam}$ is independent of the choice of $c_j$ ($j=0\nek r-1$). 
\defe{refinedtorsion}
The {\em refined torsion} of the pair $(C^\b,\Gam)$ is the element
\eq{refinedtor}
    \rho_{{}_\Gam} \ = \ \rho_{{}_{C^\b,\Gam}} \ := \  \phi_{C^\b}(c_{{}_\Gam}) \ \in \Det\big(\,H^\b(\pa)\,\big),
\end{equation}
where $ \phi_{C^\b}$ is the canonical map \refe{isomorphism}.
\edefe

\subsection{The odd signature operator}\Label{SS:oddsign}
Let $M$  be a smooth closed oriented manifold of odd dimension $d=2r-1$ and let $(E,\n)$ be a flat vector bundle over $M$. We denote by $\Ome^k(M,E)$
the space of smooth differential forms on $M$  of degree $k$ with values in $E$ and by
\[
    \n:\, \Ome^\b(M,E) \ \longrightarrow \Ome^{\b+1}(M,E)
\]
the covariant differential induced by the flat connection on $E$. Fix a Riemannian metric $g^M$ on $M$ and let $*:\Ome^\b(M,E)\to \Ome^{d-\b}(M,E)$
denote the Hodge $*$-operator. Define the {\em chirality operator} $\Gam= \Gam(g^M):\Ome^\b(M,E)\to \Ome^\b(M,E)$ by the formula
\eq{Gam}
    \Gam\, \ome \ := \ i^r\,(-1)^{\frac{k(k+1)}2}\,*\,\ome, \qquad \ome\in \Ome^k(M,E),
\end{equation}
with $r$ given as above by \/ $r=\frac{d+1}2$. The numerical factor in \refe{Gam} has been chosen so that $\Gam^2=1$, cf. Proposition~3.58 of
\cite{BeGeVe}.

\defe{oddsign}
The {\em odd signature operator} is the operator
\begin{equation} \Label{E:oddsignGam}
    \B\ = \  \B(\n,g^M) \ := \ \Gam\,\n \ + \ \n\,\Gam:\,\Ome^\b(M,E)\ \longrightarrow \  \Ome^\b(M,E).
\end{equation}
We denote by $\B_k$ the restriction of\/ $\B$ to the space $\Ome^{k}(M,E)$.
\edefe

\subsection{The graded determinant of the odd signature operator}\label{SS:grdet}
Note that for each $k=0\nek d$, the operator $\B^2$ maps $\Ome^k(M,E)$ into itself. Suppose $\calI$ is an interval of the form $[0,\lam],\
(\lam,\mu]$, or $(\lam,\infty)$ ($\mu>\lam\ge0$). Denote by $\Pi_{\B^2,\calI}$ the spectral projection of $\B^2$ corresponding to the set of
eigenvalues, whose absolute values lie in $\calI$. Set
\eq{OmecalI}\notag
    \Ome^\b_{\calI}(M,E)\ := \ \Pi_{\B^2,\calI}\big(\, \Ome^\b(M,E)\,\big)\ \subset\  \Ome^\b(M,E).
\end{equation}
If the interval $\calI$ is bounded, then, cf. Section~6.10 of \cite{BrKappelerRATdetline}, the space $\Ome^\b_\calI(M,E)$ is finite dimensional.

For each $k=0\nek d$, set
\eq{ome+-}
 \begin{aligned}
  \Ome^k_{+,\calI}(M,E) \ &:= \ \Ker\,(\n\,\Gam)\,\cap\,\Ome^k_{\calI}(M,E) \ = \ \big(\,\Gam\,(\Ker\,\n)\,\big)\,\cap\,\Ome^k_{\calI}(M,E);\\
  \Ome^k_{-,\calI}(M,E) \ &:= \ \Ker\,(\Gam\,\n)\,\cap\,\Ome^k_{\calI}(M,E)
  \ = \ \Ker\,\n\,\cap\,\Ome^k_{\calI}(M,E).
 \end{aligned}
\end{equation}
Then
\eq{Ome>0=directsum}
    \Ome^k_{\calI}(M,E) \ = \ \Ome^k_{+,\calI}(M,E) \,\oplus \Ome^k_{-,\calI}(M,E) \qquad\text{if}\quad 0\not\in \calI.
\end{equation}
We consider the decomposition \refe{Ome>0=directsum} as a {\em grading} \footnote{Note, that our grading is opposite to the one considered in
\cite[\S2]{BFK3}.} of the space $\Ome^\b_{\calI}(M,E)$, and refer to $\Ome^k_{+,\calI}(M,E)$ and $\Ome^k_{-,\calI}(M,E)$ as the positive and negative
subspaces of\/ $\Ome^k_{\calI}(M,E)$.

Set
\[
    \Ome^{\even}_{\pm,\calI}(M,E)\ =\ \bigoplus_{p=0}^{r-1}\, \Ome^{2p}_{\pm,\calI}(M,E)
\]
and let $\B^\calI$ and $\B^\calI_\even$ denote the restrictions of $\B$ to the subspaces $\Ome^\b_{\calI}(M,E)$ and $\Ome^\even_{\calI}(M,E)$
respectively. Then $\B^\calI_{\even}$ maps $\Ome^{\even}_{\pm,\calI}(M,E)$ to itself. Let $\B_{{\even}}^{\pm,\calI}$ denote the restriction of
$\B^\calI_{{\even}}$ to the space $\Ome^{\even}_{\pm,\calI}(M,E)$. Clearly, the operators $\B_{{\even}}^{\pm,\calI}$ are bijective whenever $0\not\in
\calI$.

\defe{grdetB}
Suppose $0\not\in \calI$. The {\em graded determinant} of the operator $\B^{\calI}_{\even}$ is defined by
\eq{grdetB}
    \Detgrtetnp(\B^{\calI}_{\even}) \ := \ \frac{ \Det_\tet(\B^{+,\calI}_{\even})}{\Det_\tet(-\B^{-,\calI}_{\even})}
    \ \ \in\ \ \CC\backslash \{0\},
\end{equation}
where $\Det_\tet$ denotes the $\zet$-regularized determinant associated to the Agmon angle $\tet\in (-\pi,0)$, cf., for example, \S6 of\/
\cite{BrKappelerRATdetline}.
\edefe
It follows from formula (6.17) of \cite{BrKappelerRATdetline} that \refe{grdetB} is independent of the choice of $\tet\in (-\pi,0)$.

\subsection{The canonical element of the determinant line}\Label{SS:rho}
Since the covariant differentiation $\n$ commutes with $\B$, the subspace $\Ome^\b_\calI(M,E)$ is a subcomplex of the twisted de Rham complex
$(\Ome^\b(M,E),\n)$. Clearly, for each $\lam\ge0$, the complex $\Ome^\b_{(\lam,\infty)}(M,E)$ is acyclic. Since
\eq{Ome=Ome0+Ome>0}
    \Ome^\b(M,E) \ = \  \Ome^\b_{[0,\lam]}(M,E)\,\oplus\,\Ome^\b_{(\lam,\infty)}(M,E),
\end{equation}
the cohomology $H^\b_{[0,\lam]}(M,E)$ of the complex $\Ome^\b_{[0,\lam]}(M,E)$ is naturally isomorphic to the cohomology $H^\b(M,E)$. Let
$\Gam_{\hskip-1pt{}\calI}$ denote the restriction of $\Gam$ to $\Ome^\b_\calI(M,E)$. For each $\lam\ge0$, let
\eq{thoGam0lam}
    \rho_{{}_{\Gam_{\hskip-1pt{}_{[0,\lam]}}}}\ = \ \rho_{{}_{\Gam_{\hskip-1pt{}_{[0,\lam]}}}}(\n,g^M) \in \ \Det(H^\b_{[0,\lam]}(M,E))
\end{equation}
denote the refined torsion of the finite dimensional complex $(\Ome^\b_{[0,\lam]}(M,E),\n)$ corresponding to the chirality operator
$\Gam_{\hskip-1pt{}_{[0,\lam]}}$, cf. \refd{refinedtorsion}. We view $\rho_{{}_{\Gam_{\hskip-1pt{}_{[0,\lam]}}}}$ as an element of $\Det(H^\b(M,E))$
via the canonical isomorphism between $H^\b_{[0,\lam]}(M,E)$ and $H^\b(M,E)$.

It is shown in Proposition~7.8 of \cite{BrKappelerRATdetline} that the nonzero element
\eq{rho}
    \rho(\n) \ = \ \rho(\n,g^M) \ := \ \Detgrtetnp(\B^{(\lam,\infty)}_{\even})\cdot \rho_{{}_{\Gam_{\hskip-1pt{}_{[0,\lam]}}}}
    \ \in  \ \Det(H^\b(M,E))
\end{equation}
is independent of the choice of $\lam\ge0$. Further, $\rho(\n)$ is independent of the choice of the Agmon angle $\tet\in (-\pi,0)$ of\/ $\B_\even$.
However, in general, $\rho(\n)$ might depend on the Riemannian metric $g^M$ (it is independent of $g^M$ if $\dim{}M\equiv 3\ (\MOD\ 4)$). The refined
analytic torsion, cf. \refd{rat}, is a slight modification of $\rho(\n)$, which is independent of $g^M$.

\subsection{The $\eta$-invariant}\Label{SS:etainv}
First, we recall the definition of the $\eta$-function of a non-self-adjoint elliptic operator $D$, cf. \cite{Gilkey84}. 
Let $D:C^\infty(M,E)\to C^\infty(M,E)$ be an elliptic differential operator of order $m\ge 1$ whose leading symbol is self-adjoint with respect to
some given Hermitian metric on $E$. Assume that $\tet$ is an Agmon angle for $D$ (cf., for example, Definition~3.3 of \cite{BrKappelerRAT}). Let
$\Pi_>$  (resp. $\Pi_<$) be the spectral projection whose image contains the span of all generalized eigenvectors of $D$ corresponding to eigenvalues
$\lam$ with $\RE\lam>0$ (resp. with $\RE\lam<0$) and whose kernel contains the span of all generalized eigenvectors of $D$ corresponding to
eigenvalues $\lam$ with $\RE\lam\le0$ (resp. with $\RE\lam\ge0$). For all complex $s$ with $\RE{s}<-d/m$, we define the $\eta$-function of $D$ by the
formula
\eq{eta}
    \eta_{\tet}(s,D) \ = \ \zet_\tet(s,\Pi_>,D) \ - \ \zet_\tet(s,\Pi_<,-D),
\end{equation}
where $\zet_\tet(s,\Pi_>,D):= \Tr(\Pi_>D^s)$ and, similarly, $\zet_\tet(s,\Pi_<,D):= \Tr(\Pi_<D^s)$. Note that, by the above definition, the purely
imaginary eigenvalues of $D$ do not contribute to $\eta_\tet(s,D)$.

It was shown by Gilkey, \cite{Gilkey84}, that $\eta_\tet(s,D)$ has a meromorphic extension to the whole complex plane $\CC$ with isolated simple
poles, and that it is regular at $0$. Moreover, the number $\eta_\tet(0,D)$ is independent of the Agmon angle $\tet$.

Since the leading symbol of $D$ is self-adjoint, the angles $\pm\pi/2$ are principal angles for $D$. Hence, there are at most finitely many
eigenvalues of $D$ on the imaginary axis. Let $m_+(D)$ (resp., $m_-(D)$) denote the number of eigenvalues of $D$, counted with their algebraic
multiplicities, on the positive (resp., negative) part of the imaginary axis. Let $m_0(D)$ denote the algebraic multiplicity of 0 as an eigenvalue of
$D$.

\defe{etainv}
The $\eta$-invariant $\eta(D)$ of\/ $D$ is defined by the formula
\eq{etainv}
    \eta(D) \ = \ \frac{\eta_\tet(0,D)+m_+(D)-m_-(D)+m_0(D)}2.
\end{equation}
\edefe
As $\eta_\tet(0,D)$ is independent of the choice of the Agmon angle $\tet$ for $D$, cf. \cite{Gilkey84}, so is $\eta(D)$.

\rem{etainv}
Note that our definition of $\eta(D)$ is slightly different from the one proposed by Gilkey in \cite{Gilkey84}. In fact, in our notation, Gilkey's
$\eta$-invariant is given by $\eta(D)+m_-(D)$. Hence, reduced modulo integers, the two definitions coincide. However, the number $e^{i\pi\eta(D)}$
will be multiplied by $(-1)^{m_-(D)}$ if we replace one definition by the other. In this sense, \refd{etainv} can be viewed as a {\em sign
refinement} of the definition given in \cite{Gilkey84}.
\erem

Let $\n$ be a flat connection on a complex vector bundle $E\to M$. Fix a Riemannian metric $g^M$ on $M$ and denote by
\eq{etan}
    \eta(\n) \ = \ \eta\big(\B_\even(\n,g^M)\big)
\end{equation}
the $\eta$-invariant of the restriction $\B_\even(\n,g^M)$ of the odd signature operator $\B(\n,g^M)$ to \linebreak$\Ome^\even(M,E)$.

\subsection{The refined analytic torsion}\label{E:rat}
Let $\eta_{\trivial}=\eta_{\trivial}(g^M)$ denote the $\eta$-invariant of the operator $\B_{\trivial}=\Gam\,d+d\,\Gam: \Ome^\b(M)\to \Ome^\b(M)$. In
other words, $\eta_{\trivial}$ is the $\eta$-invariant corresponding to the trivial line bundle $M\times\CC\to M$ over $M$.
\defe{rat}
Let $(E,\n)$ be a flat vector bundle on $M$. The {\em refined analytic torsion} is the element
\eq{metricindep}
    \rat\ =\ \rat(\n) \ := \ \rho(\n,g^M)\cdot \exp\,\Big(\,{i\pi\cdot\rank E}\cdot\eta_{\trivial}(g^M)\,\Big)\ \in \ \Det(H^\b(M,E)),
\end{equation}
where $g^M$ is any Riemannian metric on $M$ and $\rho(\n,g^M)\in \Det(H^\b(M,E))$ is defined by \refe{rho}.
\edefe
It is shown in Theorem~9.6 of \cite{BrKappelerRATdetline} that $\rat(\n)$ is independent of $g^M$.

\rem{rat}
In \cite{BrKappelerRAT,BrKappelerRATdetline,BrKappelerRATdetline_hol} we introduced an alternative version of the refined analytic torsion.
Consider an oriented manifold $N$ whose oriented boundary is the disjoint union of two copies of $M$. Instead of the exponential factor in
\refe{metricindep} we used the term
\[
    \exp\,\Big(\,\frac{i\pi\cdot\rank E}2\int_N\,L(p,g^M)\,\Big),
\]
where $L(p,g^M)$ is the Hirzebruch $L$-polynomial in the Pontrjagin forms of any Riemannian metric on $N$ which near $M$ is the product of\/
$g^M$ and the standard metric on the half-line. The advantage of this definition is that the latter factor is simpler to calculate than
$e^{i\pi{}\eta_{\trivial}}$. In addition, if $\dim{}M\equiv 3\ (\MOD\ 4)$, then $\int_ML(p,g^M)= 0$ and, hence, the refined analytic torsion
then coincides with $\rho(\n,g^M)$. However, in general, this version of the refined analytic torsion depends on the choice of $N$ (though only
up to a multiplication by $i^{k\cdot{\rank(E)}}$ $(k\in \ZZ)$). For this paper, however, the definition \refe{metricindep} of the refined
analytic torsion is slightly more convenient.
\erem

\subsection{Relationship with the $\eta$-invariant}\Label{SS:rat-eta}
To simplify the notation set
\eq{xi}
    T_\lam \ = \ T_\lam(\n,g^M,\tet) \ = \ \prod_{j=0}^d\,\Big(\,\Det_{2\tet} \Big[\,\big(\,
    (\Gam\n)^2+(\n\Gam)^2\,\big){\big|_{\Ome^j_{{(\lam,\infty)}}(M,E)}}\,\Big]\,\Big)^{(-1)^{j+1}\,j}.
\end{equation}
where $\tet\in (-\pi/2,0)$ and both, $\tet$ and $\tet+\pi$, are Agmon angles for $\B_\even$ (hence, $2\tet$ is an Agmon angle for $\B_\even^2$). We
shall use the following proposition, cf. \cite[Prop.~8.1]{BrKappelerRATdetline}:
\prop{xietad}
Let $\n$ be a flat connection on a vector bundle $E$ over a closed Riemannian manifold $(M,g^M)$ of odd dimension $d=2r-1$. Assume $\tet\in
(-\pi/2,0)$ is such that both $\tet$ and $\tet+\pi$ are Agmon angles for the odd signature operator $\B= \B(\n,g^M)$. Then, for every $\lam\ge0$,
\eq{xietad}
    \Big(\,\Det_{\gr,2\tet}(\B_\even^{(\lam,\infty)})\,\Big)^2 \ =\ T_\lam \cdot e^{-2\pi i\,\eta(\n,g^M)}.
\end{equation}
\eprop
Note that Proposition~8.1 of \cite{BrKappelerRATdetline}  gives a similar formula for the logarithm of $\Det_{\gr,2\tet}(\B_\even^{(\lam,\infty)})$,
thus providing a sign refined version of \refe{xietad}. In the present paper we won't need this refinement.
\prf
Set
\eq{eta=eta}
    \eta_\lam \ = \ \eta_\lam(\n,g^M) \ := \  \eta(\B_{\even}^{(\lam,\infty)}).
\end{equation}
From Proposition~8.1 and equality (10.20) of \cite{BrKappelerRATdetline} we obtain
\eq{xietad2}
    \Det_{\gr,2\tet}(\B_\even^{(\lam,\infty)})^2 \ =\ T_\lam \cdot e^{-2\pi i\,\eta_\lam}\cdot e^{-i\pi\,\dim \Ome^\even_{[0,\lam]}(M,E)}.
\end{equation}
The operator $\B^{[0,\lam]}_{\even}$ acts on the finite dimensional vector space $\Ome^\even_{[0,\lam]}(M,E)$. Hence,
$2\eta\big(\B^{[0,\lam]}_{\even}\big)\in \ZZ$ and
\eq{eta-eta}
    2\eta\big(\B^{[0,\lam]}_{\even}\big)\ \equiv \ \dim\Ome^\even_{[0,\lam]}(M,E) \quad \MOD\ 2.
\end{equation}
Since $\eta_\lam= \eta\big(\B_{\even}\big)- \eta\big(\B^{[0,\lam]}_{\even}\big)$, we obtain from \refe{eta-eta} that
\[
    e^{-i\pi\,\big(\,2\eta_\lam+\dim\Ome^\b_{[0,\lam]}(M,E)\,\big)} \ = \ e^{-2i\pi\,\eta\big(\B_{\even}\big)}.
\]
The equality \refe{xietad} follows now from \refe{xietad2}.
\eprf

\section{The Burghelea-Haller Quadratic Form}\label{SS:BH}

In this section we recall the construction of the quadratic form on the determinant line \linebreak$\Det\big(H^\b(M,E)\big)$ due to Burghelea and
Haller, \cite{BurgheleaHaller_function2}. Throughout the section we assume that the vector bundle $E\to M$ admits a non-degenerate symmetric bilinear
form $b$. Such a form, required for the construction of $\tau$, might not exist on $E$, but there always exists an integer $N$ such that on the
direct sum $E^N= E\oplus\cdots\oplus E$ of $N$ copies of $E$ such a form exists, cf. Remark~4.6 of \cite{BurgheleaHaller_function2}.

\subsection{A quadratic form on the determinant line of the cohomology of a finite dimensional complex}\label{SS:BH-fdim}
Consider the complex \refe{Cpintrod} and assume that each vector space $C^j$ ($j=0\nek d$) is endowed with a non-degenerate symmetric bilinear form
$b_j:C^j\times{}C^j \to \CC$. Set $b=\oplus{}b_j$. Then $b_j$ induces a bilinear form on the determinant line $\Det(C^j)$ and, hence, one obtains a
bilinear form on the determinant line $\Det(C^\b)$. Using the isomorphism \refe{isomorphism} we thus obtain a bilinear form on $\Det(H^\b(\pa))$.
This bilinear form induces a quadratic form on $\Det(H^\b(\pa))$, which we denote by $\tau_{C^\b,b}$.

The following lemma establishes a relationship between $\tau_{C^\b,b}$ and the construction of \refss{rtfd} and is an immediate consequence of the
definitions.
\lem{BH-rhofdim}
Suppose that $d$ is odd and that the complex $(C^\b,\pa)$ is endowed with a chirality operator $\Gam$, cf. \refss{rtfd}. Assume further that $\Gam$
preserves the bilinear form $b$, i.e. $b(\Gam x,\Gam y)= b(x,y)$, for all $x,y\in C^\b$. Then
\eq{BH-rhofdim}
    \tau_{C^\b,b}(\rho_{{}_\Gam}) \ = \ 1.
\end{equation}
where $\rho_\Gam$ is given by \refe{refinedtor}.
\elem

\subsection{Determinant of the generalized Laplacian}\label{SS:detLapl}
Assume now that $M$ is a compact oriented manifold and $E$ is a flat vector bundle over $M$ endowed with a non-degenerate symmetric bilinear form
$b$. Then $b$ together with the Riemannian metric $g^M$ on $M$ define a bilinear form
\eq{bet}
    \grb:\, \Ome^\b(M,E)\,\times \Ome^\b(M,E) \ \to \ \CC
\end{equation}
in a natural way.

Let $\n:\Ome^\b(M,E)\to \Ome^{\b+1}(M,E)$ denote the flat connection on $E$ and let $\n^\#:\Ome^\b(M,E)\to \Ome^{\b-1}(M,E)$ denote the formal
transpose of $\n$ with respect to $\grb$. Following Burghelea and Haller we define a (generalized) Laplacian
\eq{Lapl}
    \Del \ = \ \Del_{g^M,b} \ := \ \n^\#\,\n \ + \ \n\,\n^\#.
\end{equation}
Given a Hermitian metric on $E$,\ $\Del$ is not self-adjoint, but has a self-adjoint positive definite leading symbol, which is the same as the
leading symbol of the usual Laplacian. In particular, $\Del$ has a discrete spectrum, cf. \cite[\S4]{BurgheleaHaller_function2}.

Suppose $\calI$\/ is an interval of the form $[0,\lam]$ or $(\lam,\infty)$ and let $\Pi_{\Del_k,\calI}$ be the spectral projection of $\Del$
corresponding to $\calI$. Set
\eq{hatOmecalI}\notag
    \hatOme^k_{\calI}(M,E)\ := \ \Pi_{\Del_k,\calI}\big(\, \Ome^k(M,E)\,\big)\ \subset\  \hatOme^k(M,E), \qquad k=0\nek d.
\end{equation}

For each $\lam\ge0$, the  space $\hatOme^\b_{[0,\lam]}(M,E)$ is a finite dimensional subcomplex of the de Rham complex $(\Ome^\b(M,E),\n)$, whose
cohomology is isomorphic to $H^\b(M,E)$. Thus, according to \refss{BH-fdim}, the bilinear form \refe{bet} restricted to $\hatOme^\b(M,E)$ defines a
quadratic form on the determinant line $\Det(H^\b(M,E))$, which we denote by $\tau_{[0,\lam]}= \tau_{b,\n,[0,\lam]}$.

Let $\Del^\calI_k$ denote the restriction of $\Del_k$ to $\hatOme^k_\calI(M,E)$. Since the leading symbol of $\Del$ is positive definite the
$\zet$-regularized determinant $\Det'_\tet(\Del^\calI_k)$ does not depend on the choice of the Agmon angle $\tet$. Set
\eq{taulaminfty}
    \tau_{b,\n,(\lam,\infty)} \ := \ \prod_{j=0}^d\, \big(\, \Det'_\tet(\Del^{(\lam,\infty)}_j)\,\big)^{(-1)^jj} \ \in\ \CC\backslash\{0\}.
\end{equation}
Note that both, $\tau_{b,\n,[0,\lam]}$ and $\tau_{b,\n,(\lam,\infty)}$, depend on the choice of the Riemannian metric $g^M$.

\defe{BH}
The Burghelea-Haller quadratic form $\tauBH_{b,\n}$ on $\Det\big(H^\b(M,E)\big)$ is defined by the formula
\eq{BH}
    \tauBH \ = \ \tauBH_{b,\n} \ := \ \tau_{b,\n,[0,\lam]}\cdot \tau_{b,\n,(\lam,\infty)}.
\end{equation}
\edefe
It is easy to see, cf. \cite[Prop.~4.7]{BurgheleaHaller_function2}, that \refe{BH} is independent of the choice of $\lam\ge0$. Theorem~4.2 of
\cite{BurgheleaHaller_function2} states that $\tauBH$ is independent of $g^M$ and locally constant in $b$. Since we are not going to use this result in the proof of
\reft{RAT-BH}, the latter theorem provides a new proof of Theorem~4.2 of \cite{BurgheleaHaller_function2} in the case when the dimension of $M$ is odd, cf.
\refss{BHanomaly}.


\section{Proof of the Comparison Theorem}\label{S:prRAT-BH}

In this section we prove \reft{RAT-BH} adopting the arguments which we used in Section~11 of \cite{BrKappelerRATdetline} to compute the Ray-Singer
norm of the refined analytic torsion.

\subsection{The dual connection}\label{SS:dualconnection}
Suppose $M$ is a closed oriented manifold of odd dimension $d=2r-1$. Let $E\to M$ be a complex vector bundle over $M$ and let $\n$ be a flat
connection on $E$. Assume that there exists a non-degenerate bilinear form $b$ on $E$.  The {\em dual connection} $\n'$ to $\n$ with respect to the
form $b$ is defined by the formula
\[
    d\,b(u,v) \ = \ b(\n u,v) \ + \ b(u,\n' v), \qquad u,v\in C^\infty(M,E).
\]
We denote by $E'$ the  flat vector bundle $(E,\n')$.

\subsection{Choices of the metric and the spectral cut}\Label{SS:choices}
Till the end of this section we fix a Riemannian metric $g^M$ on $M$ and set $\B= \B(\n,g^M)$ and $\B'= \B(\n',g^M)$. We also fix $\tet\in
(-\pi/2,0)$ such that both $\tet$ and $\tet+\pi$ are Agmon angles for the odd signature operator $\B$. Recall that for an operator $A$ we denote by
$A^\#$ its formal transpose with respect to the bilinear form \refe{bet} defined by $g^M$ and $b$. One easily checks that
\eq{n*=2}
    \n^\# \ = \ \Gam\,\n'\,\Gam, \qquad
    (\n')^\# \ = \ \Gam\,\n\,\Gam,\qquad  \text{and}\quad \B^\# \ = \ \B',
\end{equation}
cf. the proof of similar statements when $b$ is replaced by a Hermitian form in Section~10.4 of \cite{BrKappelerRATdetline}. As $\B$ and $\B^\#$ have
the same spectrum it then follows that
\eq{n'=n}
    \eta(\B') \ = \ \eta(\B) \qquad\text{and}\qquad \Det_{\gr,\tet}(\B') \ = \ \Det_{\gr,\tet}(\B).
\end{equation}

\subsection{The duality theorem for the refined analytic torsion}\label{SS:dual}
The pairing \refe{bet} induces a non-degenerate bilinear form
\eq{Poincarepair}\notag
    H^j(M,E')\otimes H^{d-j}(M,E) \ \longrightarrow \ \CC, \qquad j=0\nek d,
\end{equation}
and, hence, identifies $H^j(M,E')$ with the dual space of $H^{d-j}(M,E)$. Using the construction of Subsection~3.4 of \cite{BrKappelerRATdetline}
(with $\tau:\CC\to \CC$ being the identity map) we thus obtain a linear isomorphism
\eq{alpE-E'}
    \alp:\, \Det\big(H^\b(M,E)\big) \ \longrightarrow \ \Det\big(H^\b(M,E')\big).
\end{equation}

We have the following analogue of Theorem~10.3 from \cite{BrKappelerRATdetline}
\th{dualityRAT}
Let\/ $E\to M$ be a complex vector bundle over a closed oriented odd-dimensional manifold $M$ endowed with a non-degenerate bilinear form $b$ and
let\/ $\n$ be a flat connection on $E$. Let\/ $\n'$ denote the connection dual to $\n$ with respect to $b$. Then
\eq{dualityRAT}
    \alp\big(\,\rat(\n)\,\big) \ = \ \rat(\n').
\end{equation}
\eth
The proof is the same as the proof of Theorem~10.3 from \cite{BrKappelerRATdetline} (actually, it is simple, since $\B$ and $\B'$ have the same
spectrum and, hence, there is no complex conjugation involved) and will be omitted.

\subsection{The Burghelea-Haller quadratic form and the dual connection}\Label{SS:RSnormdual}
Let
\eq{Laplacianprime}\notag
    \Del'\ = \ (\n')^\#\,\n' \ + \ \n'\,(\n')^\#,
\end{equation}
denote the Laplacian of the connection $\n'$. From \refe{n*=2} we conclude that
\eq{Laplace-Laplacianprime}\notag
    \Del'\ = \ \Gam\circ\Del\circ\Gam.
\end{equation}
Hence, a verbatim repetition of the arguments in Subsection~11.6 of \cite{BrKappelerRATdetline} implies that we have
\eq{dualBHlam}
    \tau_{b,\n,(\lam,\infty)} \ = \ \tau_{b,\n',(\lam,\infty)},
\end{equation}
and, for each $h\in \Det\big(H^\b(M,E)\big)$,
\eq{dualBH}
    \tauBH_{b,\n}(h) \ = \ \tauBH_{b,\n'}\big(\,\alp(h)\,\big)
\end{equation}
with $\alp$ being the duality isomorphism \refe{alpE-E'}.

From \refe{dualityRAT} and \refe{dualBH} we get
\eq{BHrho-BHrho'}
    \tauBH_{b,\n}\big(\,\rat(\n)\,\big) \ = \ \tauBH_{b,\n'}\big(\,\rat(\n')\,\big).
\end{equation}

\subsection{Direct sum of a connection and its dual}\Label{SS:n+n'}
Let
\eq{tiln}
   \tiln \ = \ \begin{pmatrix}
                        \n&0\\
                        0&\n'
                       \end{pmatrix}
\end{equation}
denote the flat connection on $E\oplus{E}$ obtained as a direct sum of the connections $\n$ and $\n'$. The bilinear form $b$ induces a bilinear form
$b\oplus{}b$ on $E\oplus{}E$. To simplify the notations we shall denote this form by $b$. For each $\lam\ge0$, one easily checks, cf. Subsection~11.7
of \cite{BrKappelerRATdetline}, that
\eq{normrattil0}
    \tau_{b,\tiln,(\lam,\infty)} \ = \ \tau_{b,\n,(\lam,\infty)}\cdot \tau_{b,\n',(\lam,\infty)}
\end{equation}
and
\eq{normrattil}
    \tauBH_{b,\tiln}\big(\,\rat(\tiln)\,\big) \ = \ \tauBH_{b,\n}\big(\,\rat(\n)\,\big) \cdot \tauBH_{b,\n'}\big(\,\rat(\n')\,\big).
\end{equation}
Combining the later equality with \refe{BHrho-BHrho'}, we get
\eq{normrattil2}
    \tauBH_{b,\tiln}\big(\,\rat(\tiln)\,\big) \ = \ \tauBH_{b,\n}\big(\,\rat(\n)\,\big)^2.
\end{equation}
Hence, {\em \refe{RAT-BH} is equivalent to the equality}
\eq{nomrattil=1}
    \tauBH_{b,\tiln}\big(\,\rat(\tiln)\,\big) \ = \ e^{-\ 4\pi i\,\big(\,\eta(\n)-\rank E\cdot\eta_{\trivial}\,\big)}.
\end{equation}

\subsection{Deformation of the chirality operator}\Label{SS:deformationGam}
We will prove \refe{nomrattil=1} by a deformation argument. For $t\in [-\pi/2,\pi/2]$ introduce the rotation $U_t$ on
\[
    \Ome^\b \ := \ \Ome^\b(M,E)\,\oplus\,\Ome^\b(M,E),
\]
given by
\[
    U_t \ = \ \begin{pmatrix}
    \ \cos t&-\,\sin t \ \\
    \ \sin t& \ \ \ \cos t \
    \end{pmatrix}.
\]
Note that $U_t^{-1}= U_{-t}$. Denote by $\tilGam(t)$ the deformation of the chirality operator, defined by
\eq{Gam(t)}
        \tilGam(t) \ = \ U_t\circ
                    \begin{pmatrix}
                        \Gam&0\\ 0&-\Gam
                    \end{pmatrix}
                    \circ U_t^{-1} \ = \ \Gam\circ
                    \begin{pmatrix}
                        \cos 2t&\sin 2t\\ \sin 2t&-\cos2t
                    \end{pmatrix}.
\end{equation}
Then
\eq{Gam0pi4}
    \tilGam(0) \ = \ \begin{pmatrix}
                        \Gam&0\\ 0&-\Gam
                    \end{pmatrix}, \qquad
    \tilGam(\pi/4) \ = \ \begin{pmatrix}
                        0&\Gam\\ \Gam&0
                    \end{pmatrix}.
\end{equation}

\subsection{Deformation of the odd signature operator}\Label{SS:deformationB}
Consider a one-parameter family of operators $\tB(t):\Ome^\b\to \Ome^\b$ with $t\in [-\pi/2,\pi/2]$ defined by the formula
\eq{tB(t)}
    \tB(t) \ := \ \tilGam(t)\,\tiln\ + \ \tiln\,\tilGam(t).
\end{equation}
Then
\eq{tilDel0}
    \tB(0) \ = \ \begin{pmatrix}
                        \B&0\\ 0&-\B'
                    \end{pmatrix}
\end{equation}
and
\eq{tilBpi/4}
    \tB(\pi/4) \ = \ \begin{pmatrix}
                        0&\Gam\n'+\n\Gam\\ \Gam\n+\n'\Gam&0
                    \end{pmatrix}.
\end{equation}
Hence, using \refe{n*=2}, we obtain
\eq{tilDelpi}
                    \tB(\pi/4)^2 \ = \ \begin{pmatrix}
                        \Del&0\\ 0&\Del'
                    \end{pmatrix} \ = \ \tilDel.
\end{equation}
Set
\begin{gather}
    \Ome^\b_+(t) \ := \ \Ker\, \tiln\,\tilGam(t);\notag\\
        \Ome^\b_- \ := \ \Ker\tiln \ = \ \Ker \n\oplus \Ker\n'.\notag
\end{gather}
Note that $\Ome^\b_-$ is independent of $t$. Since the operators $\tiln$ and $\tilGam(t)$ commute with $\tB(t)$, the spaces $\Ome^\b_+(t)$ and
$\Ome^\b_-$ are invariant for $\tB(t)$.

Let $\calI$ be an interval of the form $[0,\lam]$ or $(\lam,\infty)$. Denote
\eq{Omecal(t)I}\notag
    \Ome^\b_{\calI}(t)\ := \ \Pi_{\tB(t)^2,\calI}\big(\, \Ome^\b(t)\,\big)\ \subset\  \Ome^\b(t),
\end{equation}
where $\Pi_{\tB(t)^2,\calI}$ is the spectral projection of $\tB(t)^2$ corresponding to $\calI$. For $j=0\nek d$, set $\Ome^j_\calI(t)=
\Ome^\b_\calI(t)\cap\Ome^j$ and
\eq{Ome+-tlam}
    \Ome^j_{\pm,\calI}(t) \ := \ \Ome^j_\pm(t)\cap \Ome^j_\calI(t).
\end{equation}

As $\Pi_{\tB(t)^2,\calI}$ and $\tB(t)$ commute, one easily sees, cf. Subsection~11.9 of \cite{BrKappelerRATdetline}, that
\eq{Omet+-}
    \Ome^\b_{(\lam,\infty)}(t) \ = \ \Ome^\b_{+,(\lam,\infty)}(t)\, \oplus\, \Ome^\b_{-,(\lam,\infty)}(t),\qquad t\in [-\pi/2,\pi/2].
\end{equation}

We define $\tB_j^\calI(t),\ \tB_\even^\calI(t),\ \tB_\odd^\calI(t),\ \tB_j^{\pm,\calI}(t),\ \tB_\even^{\pm,\calI}(t),\ \tB_\odd^{\pm,\calI}(t)$, etc.
in the same way as the corresponding maps were defined in \refss{grdet}.

\subsection{Deformation of the canonical element of the determinant line}\Label{SS:deformationrho}
Since the operators $\tiln$ and $\tB(t)^2$ commute, the space $\Ome^\b_{\calI}(t)$ is invariant under $\tiln$, i.e., it is a subcomplex of $\Ome^\b$.
The complex $\Ome^\b_{(\lam,\infty)}(t)$ is acyclic and, hence, the cohomology of the finite dimensional complex $\Ome^\b_{[0,\lam]}(t)$ is naturally
isomorphic to
\[
    H^\b(M,E\oplus{}E')\ \simeq\ H^\b(M,E)\oplus{}H^\b(M,E').
\]
Let $\tilGam_{[0,\lam]}(t)$ denote the restriction of $\tilGam(t)$ to $\Ome^\b_{[0,\lam]}(t)$. As $\tilGam(t)$ and $\tB(t)^2$ commute, it follows
that $\tilGam_{[0,\lam]}(t)$ maps $\Ome^\b_{[0,\lam]}(t)$ onto itself and, therefore, is a chirality operator for $\Ome^\b_{[0,\lam]}(t)$. Let
\eq{rhotilGam}
    \rho_{{}_{\tilGam_{\hskip-1pt[0,\lam]}(t)}}(t) \ \in \ \Det\big(H^\b(M,E\oplus E')\,\big)
\end{equation}
denote the refined torsion of the finite dimensional complex $\big(\Ome^\b_{[0,\lam]}(t),\tiln\big)$ corresponding to the chirality operator
$\tilGam_{[0,\lam]}(t)$, cf. \refd{refinedtorsion}.

{}For each $t\in (-\pi/2,\pi/2)$ fix an Agmon angle $\tet= \tet(t)\in (-\pi/2,0)$ for $\tB_{\even}(t)$  and define the element $\rho(t)\in
\Det\big(\,H^\b(M,E\oplus E')\,\big)$ by the formula
\eq{tilrho(t)}
    \rho(t) \ := \ \Detgrtetnp\big(\tB^{(\lam,\infty)}_{\even}(t)\big)\cdot \rho_{{}_{\tilGam_{\hskip-1pt{}_{[0,\lam]}}(t)}}(t),
\end{equation}
where $\lam$ is any non-negative real number. It follows from Proposition~5.10 of \cite{BrKappelerRATdetline} that $\rho(t)$ is independent of the
choice of $\lam\ge0$.

For $t\in [-\pi/2,\pi/2]$, $\lam\ge0$, set
\eq{xilam(t)}
    T_\lam(t) \ := \ \prod_{j=0}^d\,
    \Big(\, \Det_{2\tet}\big[\,{\tB^{(\lam,\infty)}_{\even}(t)^2}{\big|_{\Ome^j_{{(\lam,\infty)}}(t)}}\,\big]\,\Big)^{(-1)^{j+1}\,j},
\end{equation}
Then, from \refe{tilrho(t)} and \refe{xietad} we conclude that
\eq{normrho(t)=eta}
    \tauBH_{b,\tiln}\big(\,\rho(t)\,\big) \ = \ \tauBH_{b,\tiln}\Big(\,\rho_{{}_{\tilGam_{\hskip-1pt{}_{[0,\lam]}}(t)}}(t)\,\Big)\cdot T_\lam(t)
        \cdot
    e^{-2i\pi\,\eta\big(\tB_{\even}(t)\big)}.
\end{equation}
In particular,
\[
    \tauBH_{b,\tiln}\Big(\,
        \rho_{{}_{\tilGam_{\hskip-1pt{}_{[0,\lam]}}(t)}}(t)\,\Big)\cdot T_\lam(t)
\]
is independent of $\lam\ge0$.

\subsection{Computation for $t=0$}\Label{SS:rho(0)}
{}From \refe{Gamc} and definition \refe{refinedtor} of the element $\rho$, we conclude that
\eq{rho-Gam}\notag
    \rho_{{}_{-\Gam_{\hskip-1pt{}_{[0,\lam]}}}}(\n',g^M) \ = \ \pm\,\rho_{{}_{\Gam_{\hskip-1pt{}_{[0,\lam]}}}}(\n',g^M).
\end{equation}
Thus,
\eq{tautil-tautil'}\notag
    \tauBH_{b,\tiln}\big(\,\rho_{{}_{-\Gam_{\hskip-1pt{}_{[0,\lam]}}}}(\n',g^M)\,\big) \ = \
    \tauBH_{b,\tiln}\big(\,\rho_{{}_{\Gam_{\hskip-1pt{}_{[0,\lam]}}}}(\n',g^M)\,\big).
\end{equation}
Hence, from \refe{tiln} and \refe{Gam0pi4} we obtain
\eq{tau(r(o))}
    \tauBH_{b,\tiln}\big(\,\rho_{{}_{\tilGam_{\hskip-1pt{}_{[0,\lam]}}(0)}}(0)\,\big) \ = \
    \tauBH_{b,\n}\big(\,\rho_{{}_{\Gam_{\hskip-1pt{}_{[0,\lam]}}}}(\n,g^M) \, \big) \cdot
    \tauBH_{b,\n'}\big(\,\rho_{{}_{\Gam_{\hskip-1pt{}_{[0,\lam]}}}}(\n',g^M)\,\big).
\end{equation}
Using \refe{tilDel0} and the definitions \refe{xi} and \refe{xilam(t)} of $T_\lam$ we get
\eq{Tlam-Tlam'}
    T_\lam(0) \ = \ T_\lam(\n,g^M,\tet)\cdot T_\lam(\n',g^M,\tet).
\end{equation}

Combining the last two equalities with definitions \refe{rho}, \refe{tilrho(t)} of $\rho$ and with \refe{xietad}, \refe{n'=n}, and
\refe{BHrho-BHrho'}, we obtain
\eq{rho(0)'}
    \tauBH_{b,\tiln}\big(\,
        \rho_{{}_{\tilGam_{\hskip-1pt{}_{[0,\lam]}}(0)}}(0)\,\big)\cdot T_\lam(0) \ = \
        \tauBH_{b,\n}\big(\,\rat(\n)\,\big)^2\cdot e^{ 4\pi i\,\big(\,\eta(\n)-\rank E\cdot\eta_{\trivial}\,\big)}.
\end{equation}
Comparing this equality with \refe{normrattil2} we see that {\em in order to prove \refe{nomrattil=1} and, hence, \refe{RAT-BH} it is enough to show
that}
\eq{rho(0)'final}
    \tauBH_{b,\tiln}\big(\,
        \rho_{{}_{\tilGam_{\hskip-1pt{}_{[0,\lam]}}(0)}}(0)\,\big)\cdot T_\lam(0) \ = \ 1.
\end{equation}

\subsection{Computation for $t=\pi/4$}\Label{SS:rho(pi/4)}
From \refe{tilDelpi} and the definitions \refe{taulaminfty} and \refe{xilam(t)} of $\tau_{b,\tiln,(\lam,\infty)}$ and $T_\lam(t)$, we conclude
\eq{tB(pi/4)=TRS}
    T_\lam(\pi/4) \ = \ 1/\tau_{b,\tiln,(\lam,\infty)}.
\end{equation}

By \refe{tilDelpi} we have
\[
    \Ome^\b_{[0,\lam]}(\pi/4) \ = \ \Ome^\b_{[0,\lam]}(M,E)\oplus \Ome^\b_{[0,\lam]}(M,E').
\]
From \refe{Gam0pi4} we see that the restriction of $\tilGam(\pi/4)$ to $\Ome^\b_{[0,\lam]}(\pi/4)$ preserves the bilinear form on
$\Ome^\b_{[0,\lam]}(\pi/4)$ induced by $b$. Hence we obtain from \refl{BH-rhofdim}
\eq{normrholam(pi)}\notag
    \tau_{b,\tiln,[0,\lam]}\big(\,\rho_{{}_{\tilGam_{\hskip-1pt[0,\lam]}(\pi/4)}}(\pi/4)\,\big) \ = \ 1.
\end{equation}
Therefore, from \refe{tB(pi/4)=TRS} and the definitions \refe{BH} of $\tauBH$, we get
\eq{normrho(pi/4)}
     \tauBH_{b,\tiln}\big(\,
        \rho_{{}_{\tilGam_{\hskip-1pt{}_{[0,\lam]}}(\pi/4)}}(\pi/4)\,\big)\cdot T_\lam(\pi/4) \ = \ 1.
\end{equation}

\subsection{Proof of \reft{RAT-BH}}\Label{SS:prRAT-BH}
Fix an Agmon angle $\tet\in (-\pi/2,0)$ and set
\[
    \xi_{\lam,\tet}(t) \ := \ -\,\frac12\,\sum_{j=0}^d\, (-1)^{j+1}\,j\, \zet_\tet'\big(\,0,\tilB_\even(t)^2|_{\Ome_{(\lam,\infty)}^j(t)}\,\big),
\]
where $\zet_\tet'(0,A)$ denotes the derivative at zero of the $\zet$-function of the operator operator $A$. Then $T_\lam(t)=
e^{2\xi_{\lam,\tet}(t)}$. Hence, from \refe{normrho(pi/4)} we conclude that {\em in order to prove \refe{rho(0)'final} (and, hence,
\refe{nomrattil=1} and \refe{RAT-BH}) it suffices to show that
\eq{normrholam(t)xilam(t)}
     \tauBH_{b,\tiln}\big(\,
        \rho_{{}_{\tilGam_{\hskip-1pt{}_{[0,\lam]}}(t)}}(t)\,\big)\cdot e^{2\xi_{\lam,\tet}(t)}
\end{equation}
is independent of\/ $t$. }

{}Fix $t_0\in[-\pi/2,\pi/2]$ and let $\lam\ge0$ be such that the operator $\tB_{\even}(t_0)^2$ has no eigenvalues with absolute value $\lam$. Choose
an angle $\tet\in (-\pi/2,0)$ such that both $\tet$ and $\tet+\pi$ are Agmon angles for $\tB(t_0)$. Then there exists $\del>0$ such that for all
$t\in (t_0-\del,t_0+\del)\cap[-\pi/2,\pi/2]$, the operator $\tB_{\even}(t)^2$ has no eigenvalues with absolute value $\lam$ and both $\tet$ and
$\tet+\pi$ are Agmon angles for $\tB(t)$.

A verbatim repetition of the proof of Lemma~9.2 of \cite{BrKappelerRATdetline} shows that
\eq{ddttilxilam}
    \frac{d}{dt}\, \rho_{{}_{\tilGam_{\hskip-1pt{}_{[0,\lam]}}(t)}}(t)\cdot e^{\xi_{\lam,\tet}(t)} \ = \ 0.
\end{equation}
Hence, \refe{normrholam(t)xilam(t)} is independent of $t$. \hfill$\square$

\section{Properties of the Burghelea-Haller Quadratic Form}\label{SS:propBH}

Combining \reft{RAT-BH} with results of our papers \cite{BrKappelerRAT,BrKappelerRATdetline,BrKappelerRATdetline_hol} we derive new properties and
obtain new proofs of some known ones of the Burghelea-Haller quadratic form $\tau$. In particular, we prove a weak version of the Burghelea-Haller
conjecture, \cite[Conjecture~5.1]{BurgheleaHaller_function2}, which relates the quadratic form \refe{BHmodified} with the Farber-Turaev torsion.

\subsection{Independence of{} $\tauBH$ of the Riemannian metric and the bilinear form}\label{SS:BHanomaly}
The following theorem was established by Burghelea and Haller
\cite[Th.~4.2]{BurgheleaHaller_function2} without the assumption that $M$ is oriented and odd-dimensional.
\th{BHanomaly}{\rm\textbf{[Burghelea-Haller]}}
Let $M$ be an odd dimensional orientable closed manifold and let $(E,\n)$ be a flat vector bundle over $M$. Assume that there exists a non-degenerate symmetric
bilinear form $b$ on $E$. Then the Burghelea-Haller quadratic form $\tauBH_{b,\n}$ is independent of the choice of the Riemannian metric $g^M$ on $M$ and is locally
constant in $b$.
\eth
Our \reft{RAT-BH} provides a new proof of this theorem and at the same time gives the following new result.
\th{BHanomaly2}
Under the assumptions of \reft{BHanomaly} suppose that\/ $b'$ is another non-degenerate symmetric bilinear form on $E$ not necessarily homotopic to $b$ in the space
of non-degenerate symmetric bilinear forms. Then $\tauBH_{b',\n} = \pm\tauBH_{b,\n}$.
\eth
\subsubsection*{Proof of Theorems~\ref{T:BHanomaly} and \ref{T:BHanomaly2}}\label{SS:prBHanomaly} As the refined analytic torsion $\rat(\n)$ does not depend on $g^M$
and $b$, \reft{RAT-BH} implies that, modulo sign, $\tauBH_{b,\n}$ is independent of $g^M$ and $b$. Since $\tauBH_{b,\n}$ is continuous in $g^M$ and $b$ it follows
that it is locally constant in $g^M$ and $b$. Since the space of Riemannian metrics is connected, $\tauBH_{b,\n}$ is independent of $g^M$. \hfill$\square$

\subsection{Comparison with the Farber-Turaev torsion: proof of \reft{BHconj}}\label{SS:comparisonBHTUR}
Let $L(p)= L_M(p)$ denote the Hirzebruch $L$-polynomial in the Pontrjagin forms of a Riemannian metric on $M$.  We write $\hatL(p)\in H_\b(M,\ZZ)$
for the Poincar\'e dual of the cohomology class $[L(p)]$ and let $\hatL_1\in H_1(M,\ZZ)$ denote the component of $\hatL(p)$ in $H_1(M,\ZZ)$.

Theorem~5.11 of \cite{BrKappelerRATdetline_hol} combined with formulae (5.4) and (5.6) of \cite{BrKappelerRATdetline_hol} implies that for each connected component
$\C\subset \Flat(E)$, there exists a constant $F_\C$ such that for every flat connection $\n\in \C$ and every Euler structure $\eps$ we have
\eq{FoC}
    |F_\C| \ = \ \big|\,e^{-2\pi i\<\Arg_\n,\hatL_1\>\, +\,\ 2\pi i\,\eta(\n)}\,\big|,
\end{equation}
and
\eq{RAT-Tur}
    \left(\,\frac{\rho_{\eps,\gro}(\n)}{\rat(\n)}\,\right)^2 \ = \ F_\C\cdot e^{2\pi i\<\Arg_\n,\hatL_1+c(\eps)\>}.
\end{equation}
Hence, from the definition \refe{ratqf} of the quadratic form $\tau$, we  get
\eq{BH-Tur2}
     \tau_\n\big(\,\rho_{\eps,\gro}(\n)\,\big)\cdot   e^{-2\pi i\<\Arg_\n,c(\eps)\>} \ = \  F_\C\cdot
        e^{2\pi i\<\Arg_\n,\hatL_1\>\, -\,\ 2\pi i\,(\,\eta(\n)-\rank E\cdot\eta_{\trivial}\,)}.
\end{equation}

Assume now that $\n_t$ with $t\in [0,1]$ is a smooth path of flat connections. The derivative $\dot\n_t= \frac{d}{dt}\n_t$ is a smooth differential
1-form with values in the bundle of isomorphisms of $E$. We denote by $[\Tr\dot\n_t]\in H^1(M,\CC)$ the cohomology class of the closed 1-form
$\Tr\dot\n_t$.

By Lemma~12.6 of \cite{BrKappelerRAT}, we have
\eq{ddtArg}
        2\pi i \,\frac{d}{dt}\,\Arg_{\n_t} \ = \ -\big[\Tr\,\dot\n_t\big] \  \in \ H^1(M,\CC).
\end{equation}

Let $\oeta(\n_t,g^M)\in \CC/\ZZ$ denote the reduction of $\eta(\n_t,g^M)$ modulo $\ZZ$. Then $\oeta(\n_t,g^M)$ depends smoothly on $t$, cf.
\cite[\S1]{Gilkey84}. From Theorem~12.3 of \cite{BrKappelerRAT} we obtain\footnote{This result was originally proven by Gilkey
\cite[Th.~3.7]{Gilkey84}.}
\eq{ddteta}
    -\,2\pi i\,\frac{d}{dt}\,\oeta(\n_t,g^M) \ = \ \int_M\, L(p)\wedge\Tr\dot\n_t \ = \  \big\<\big[\Tr\dot\n_t\big],\hatL_1\big\>.
\end{equation}

From \refe{BH-Tur2}--\refe{ddteta} we then obtain
\eq{ddttaurho}
    \frac{d}{dt}\,\Big[\,\tau_{\n_t}\big(\,\rho_{\eps,\gro}(\n_t)\,\big)\cdot   e^{-2\pi i\<\Arg_\n,c(\eps)\>}\,\Big] \ = \ 0,
\end{equation}
proving that the right hand side of \refe{BH-Tur2} is independent of\/ $\n\in \C$. From \refe{FoC} and the fact that $\eta_{\trivial}\in \RR$ we
conclude that the absolute value of the right hand side of \refe{BH-Tur2} is equal to 1. Part (i) of \reft{BHconj} is proven.

Finally, consider the case when $\C$ contains an acyclic Hermitian connection $\n$. In this case both, $\tau_\n$ and $\rho_{\eps,\gro}(\n)$, can be viewed as non-zero
complex numbers. To prove part (ii) of \reft{BHconj} it is now enough to show that the numbers $\rho_{\eps,\gro}(\n)^2$ and $\tau_\n\cdot{}e^{-2\pi
i\<\Arg_\n,c(\eps)\>}$ have the same phase. Since $\n$ is a Hermitian connection, the number $\eta(\n)$ is real. Hence, it follows from Theorem~10.3 of
\cite{BrKappelerRATdetline} that
\[
    \Ph\big(\,\rat(\n)\,\big) \ \equiv \ -\, \pi i\,\big(\,\eta(\n)-\rank E\cdot\eta_{\trivial}\,\big) \qquad \MOD \pi i.
\]
Thus, by \refe{ratqf},
\eq{Phasetau}
    \Ph\big(\,\tau_\n\cdot{} e^{-2\pi i\<\Arg_\n,c(\eps)\>}\,\big) \ = \ \Ph\big(\, e^{-2\pi i\<\Arg_\n,c(\eps)\>}\,\big)
    \ = \ -\,2\pi\,\RE\,\<\Arg_\n,c(\eps)\>.
\end{equation}
By formula (2.4) of \cite{FarberTuraev99},
\eq{Phaserho}
    \Ph\big(\,\rho_{\eps,\gro}(\n)^2\,\big) \ = \ -\,2\pi\,\RE\,\<\Arg_\n,c(\eps)\>.
\end{equation}
The proof of \reft{BHconj} is complete. \hfill$\square$

\providecommand{\bysame}{\leavevmode\hbox to3em{\hrulefill}\thinspace} \providecommand{\MR}{\relax\ifhmode\unskip\space\fi MR }
\providecommand{\MRhref}[2]{%
  \href{http://www.ams.org/mathscinet-getitem?mr=#1}{#2}
} \providecommand{\href}[2]{#2}

\end{document}